\documentclass[pdftex,10pt,a4paper]{article}
\usepackage{times}
\usepackage[pdftex]{graphicx}

\usepackage{amssymb} 

\usepackage{amsmath} 

\usepackage{amsthm} 

\usepackage[center]{titlesec}

\theoremstyle{plain}

\newtheorem*{mainth}{Theorem}
\theoremstyle{definition}

\begin{document}



\begin{center}

\Large{MA\L GORZATA STAWISKA}\footnote{dr Ma\l gorzata Stawiska-Friedland, Mathematical Reviews, 416 Fourth St., Ann Arbor, MI 48103, USA\\email: stawiska@umich.edu}\\

\huge{\bf  LUCJAN EMIL B\"OTTCHER (1872-1937)--THE POLISH PIONEER OF HOLOMORPHIC DYNAMICS}\\
-----------------------------------\\

\huge{\bf LUCJAN EMIL B\"OTTCHER (1872-1937)--POLSKI PIONIER DYNAMIKI HOLOMORFICZNEJ}\\
\end{center}

\noindent \normalsize{{\it A b s t r a c t} In this article I will present Lucjan Emil B\"ottcher (1872-1937),  a little-known Polish mathematician active in Lw\'ow.  I will outline his scholarly path and  describe briefly his mathematical achievements. In view of later developments in holomorphic dynamics, I will argue that, despite some flaws in his work,  B\"ottcher should be regarded not only as a contributor to the area but in fact as one of its founders.}\\


\noindent \normalsize{{\it S t r e s z c z e n i e} W  artykule tym przedstawi\c e Lucjana Emila B\"ottchera (1872-1937), ma\l o znanego matematyka polskiego dzia\l aj\c acego we Lwowie. Naszkicuj\c e jego drog\c e naukow\c a    oraz  opisz\c e  pokr\'otce jego osi\c agni\c ecia matematyczne. W \'swietle p\'o\'zniejszego rozwoju dynamiki holomorficznej  b\c ed\c e\ dowodzi\'c, \.ze mimo pewnych wad jego prac  nale\.zy nie tylko doceni\'c wk\l ad B\"ottchera w t\c e dziedzin\c e, ale w istocie zaliczy\'c go  do jej tw\'orc\'ow.}\\

{\it S\l owa kluczowe: matematyka XIX i  XX wieku, dynamika holomorficzna, iteracje funkcji wymiernych}\\

 \newpage

\normalsize

\section{Introduction}

Holomorphic dynamics-- in particular  the study of iteration of rational maps on the Riemann sphere-- is an active area of current mathematical research. Among many mathematicians who have worked  in it there are several recipients of the Fields medal (the highest honor in mathematics): John Milnor, William Thurston, Curtis McMullen and Jean-Christophe Yoccoz.  While the theory requires rather deep knowledge of concepts and methods of various areas of mathematics, including complex function theory,  dynamical systems, topology, number theory, etc., its objects have become to some extent present in the popular culture, thanks to  amazing computer-generated pictures of Julia sets and the Mandelbrot set.  Some problems in holomorphic dynamics can be traced back to Arthur Cayley or even Isaac Newton, but the common view is that it got its start as a systematic and separate  area of mathematics around 1918, with works of Pierre Fatou, Gaston Julia, Samuel  Latt\'es and Salvatore Pincherle. The beginnings are discussed in two recent historical monographs: one by Mich\`ele Audin (\cite{Au}) and one by Daniel Alexander, Felice Iavernaro and Alessandro Rosa (\cite{AIR}).  These publications, as well as  recent  textbooks  on the subject of holomorphic dynamics, mention Lucjan Emil B\"ottcher among earlier contributors to the area.  His result on the local behavior of a holomorphic function near its superattracting fixed point, referred to as B\"ottcher's theorem (in which so-called B\"ottcher's equation and B\"ottcher's function appear) is a classic one:

\begin{mainth} Let $f(z) = a_mz^m+a_{m+1}z^{m+1}+..., \quad m \geq 2, a_m \neq 0$ be an analytic function in a neighborhood of $0$. Then there exists a conformal map $F$ of a neighborhood of $0$ onto the unit disk, $F(z)=z+bz^2+...$,  satisfying the equation $Ff(z) = [F(z)]^m$.
\end{mainth}

This note concerns  Lucjan Emil B\"ottcher, the  author of the above result, a Polish mathematician active at the turn of the 19th/20th centuries, who until recently remained a rather obscure figure.  Commonly available biographical information about him came from his curriculum vitae which he wrote in Latin and attached to his doctoral thesis published in 1898, so almost nothing was known about his later years. \footnote{Stanis\l aw Domoradzki discovered many materials concerning B\"ottcher and gave them scholarly treatment in his recent publications \cite{D1} and \cite{D2}.  He and I carried out an extensive analysis of B\"ottcher's work  in our joint article \cite{DS}.}  His mathematical output also became largely forgotten. One can point to several possible reasons why this could happen. First, B\"ottcher wrote mostly in Polish, so his work could not be widely read  by international mathematical community (his most cited paper was written in Russian.) Second, he was active in Lw\'ow (later Lvov, now Lviv),  an academic center frequently caught in the turbulent history of the 20th century (wars, changing political borders and governments), so many documents of his activity could not be easily retrieved; some might have been dispersed or lost. Third, a reason more interesting for a historian of mathematics, B\"ottcher's results were not appreciated by his contemporaries from academic establishment. He was an academic teacher, but never made it to the rank of a professor, so he could not disseminate his ideas by guiding doctoral students or holding specialized lectures and seminars (although, as a part of his habilitation proceedings, he proposed a plan of lectures on the general theory of iteration).\\

One should note that B\"ottcher's work was quite removed from  interests of Lw\'ow mathematicians of that time.  Worse, it  contained flaws: most proofs were only sketched, some conclusions were unjustified and the notions were not always well defined. His critics did not see a wealth of ideas, examples and partial results which amounted to  almost complete outline of the theory  developed independently only some 20 years later  by Fatou, Julia, Latt\'es and Pincherle. In what follows I will  talk about perception of B\"ottcher's work by his contemporaries. I will go back to his doctoral  studies in Leipzig with Sophus Lie and mention a controversy regarding the evaluation of B\"ottcher's thesis, in which Lie stood  for his student against his academic colleagues. I will also discuss B\"ottcher's later academic career in Lw\'ow, in particular his  repeated-- but ultimately unsuccesful-- attempts to obtain habilitation at the Lw\'ow University.  It was only after 1920 that importance of B\"ottcher's results was realized by other mathematicians (starting with Joseph Fels Ritt, an American who published the first complete proof of B\"ottcher's theorem).  Nowadays B\"ottcher's theorem is well known to  researchers in holomorphic dynamics and functional equations (and was generalized in many ways, see \cite{DS}), and he rightly gets the credit for constructing the first example of an everywhere chaotic rational map. But there is more to Lucjan Emil B\"ottcher's mathematical output that needs to be better known and appreciated.\\

\section{The life  of Lucjan Emil B\"ottcher}

Lucjan Emil B\"ottcher was born in Warsaw on January 7 (21 according to the new style calendar), 1872, in a family of Evangelical-Lutheran denomination. He attended private real schools in Warsaw and passed his maturity exam in 1893 in the (classical) gymnasium in \L om\.za. The same year he enrolled in the Division of Mathematics and Physics of the Imperial University of Warsaw (where Russian was the language of instruction),  attending lectures in mathematics, astronomy, physics and chemistry. In 1894 he was expelled from the university for participating in a Polish patriotic manifestation. He left Warsaw and moved to Lw\'ow.  He became a student   in the Division of Machine Construction at the Lw\'ow Polytechnic School,  where in 1896 he passed a state exam  obtaining a certificate with distinction, and got his so-called half-diploma in 1897. Wishing to complete a course of university studies in mathematics, he then moved to Leipzig. He spent  three semesters at the university there, attending lectures  in mathematics, physics and psychology. He completed his studies getting the degree of doctor of philosophy (under the direction of Sophus Lie, one of the most important mathematicians of the 19th century) in 1898 on the basis of the doctoral thesis ''Beitr\"age zu der Theorie der Iterationsrechnung" as well as examinations in mathematics, geometry and physics.\\

After finishing his studies B\"ottcher returned to Lw\'ow and took a job of an assistant in the (Imperial and Royal) Lw\'ow Polytechnic School (initially at the Chair of Mechanical Technology, then at the reactivated Chair of Mathematics). He had his PhD diploma nostrified in 1901. In 1910 he became an \textit{adiunkt} and in 1911 he obtained the license to lecture  (\textit{venia legendi}) and habilitation in mathematics in the Lw\'ow Polytechnic School. He lectured on many mathematical subjects in the engineering curriculum as well as on theoretical mechanics. Starting in 1901 he made several attempts to obtain habilitation also at the Lw\'ow University, all of which were unsuccesful.\\

Besides academic teaching, B\"ottcher's activities comprised taking part in meetings and conventions of mathematicians and philosophers (he was a member of the Polish  Mathematical Society founded in 1917), writitng and publishing papers in mathematics, mathematics education, logic and mechanics. He also wrote lecture notes, textbooks for the use in high schools and  booklets on spiritualism and afterlife. The (so far most) complete list of his publications can be found in \cite{DS}.

Lucjan Emil B\"ottcher retired from the Lw\'ow Polytechnic School in 1935 and died in Lw\'ow on May 29, 1937.\\

\section{B\"ottcher in Leipzig}

 Lucjan Emil B\"ottcher enrolled in the University of Leipzig on February 1, 1897, in order to study mathematics. He took courses from Sophus Lie (theory of differential invariants, theory of differential equations with known infinitesimal transforms, theory of continuous groups of transformations, seminars in theory of integral invariants and in differential equations), Adolph Mayer (higher analytic mechanics), Friedrich Engel (differential equations, algebraic equations, non-euclidean geometry), Felix Hausdorff (similarity transformations), Paul Drude (electricity and magnetism), Gustaw Wiedemann (exercises in physics) and Wilhelm Wundt (psychology).\footnote{Information about B\"ottcher's coursework is taken from the curriculum vitae which he wrote in Polish in 1901, now in the archives of Lviv University.}  On February 7, 1898, B\"ottcher submitted his dissertation ``Beitr\"age zu der Theorie der Iterationsrechnung".  His supervisor  was Sophus Lie--  the creator of the theory of continuous groups of transformations and one of the most important mathematicians of the 19th century.\footnote{Sophus Lie (1842-1899) was of Norwegian nationality. His work made an impact on the development of modern geometry, algebra and differential equations. His results led to the  emergence of new areas of mathematics, e.g., topological groups. They  remain significant for today's mathematics  and also found applications in quantum physics.}.
Like many mathematicians who came into contact with Lie (in Leipzig and elsewhere-- e.g., Victor B\"acklund, \'Elie Cartan, Wilhelm Killing, Emile Picard, Henri Poincar\'e, Eduard Study and Ernest Vessiot), B\"ottcher became fascinated with Lie's theories. He wanted to define iterations of maps with arbitrary exponents and then to study relations between iterations and functional equations, and he thought of using the theory of continuous groups of  transformations as the framework for his considerations. In Chapter I of his  dissertation he expressed some  formal relations for iterations  by means of one-parameter continuous groups of transformations. He formulated some ``fundamental theorems" without proofs, claiming only that they were special cases of some results by Lie. Part II  of B\"ottcher's  dissertation was devoted to the study of iterations of rational functions of a complex variable (ranging over the Riemann sphere) and contain results and ideas which can be regarded as foundations of holomorphic dynamics (see further  sections or \cite{DS} for more information). In part  III B\"ottcher resumed the general study of relations between iteration and functional equations. The first draft of the dissertation was ready in January 1898 and was then presented to Lie.\\

However,  on March 1 B\"ottcher still could not be admitted to  doctoral examinations. According to the handwritten note at the bottom of the document confirming the opening of the official proceedings for the degree of doctor of philosophy \footnote{The documents about B\"ottcher's doctoral proceedings come from the archives of University of Leipzig. All documents cited in this paper were  found by S. Domoradzki.}, signed by the pro-chancellor Ferdinand Zirkel,  Wilhelm Scheibner refused to submit a report on B\"ottcher's dissertation, and at the urgent request of Lie a second examiner had to be found. (On the same document it can be seen that Scheibner's name as an examiner  was struck out, and the name of Mayer was appended after the name of Lie.) This ``urgent request" was in fact made in writing. Here are  excerpts from Lie's note:\\

\textit{``At present I cannot recognize that the author has  definitely managed to substantiate significant new results. Despite all of this, his considerations, which testify to  diligence and talent, have their value. (...)\\
``In any case, I (as well as Mr. Scheibner) agree that this attempt be accepted as a thesis, and we also agree regarding the evaluation being II. I choose such a good grade because Mr. B\"ottcher himself chose his topic and developed it independently.(...)\\
``Under the conditions mentioned above, I support the acceptance of the dissertation with evaluation II and admission to the oral exam".}\\

A further document dated April 27, 1898 records B\"ottcher's completion of the required examinations in mathematics and physics, along with the examiners' evaluations and overall grade ``magna cum laude" (IIa), as well as his promotion to the degree.  Sophus Lie's written evaluation was that ``[t]he candidate is an intelligent mathematician, possessing good and solid knowledge." There was also a longer examination report by Adolph Mayer.\\

Lie  recognized the failure of B\"ottcher's initial goal of grounding the theory of iteration in the theory of continuous groups of transformations. Nevertheless, his opinion on B\"ottcher and his achievements was high. This is especially noteworthy, as out of 56 students in Leipzig who completed their doctorates in mathematics between 1890 and 1898, 26 did so with Lie (including his later  collaborator  Georg Scheffers, as well as  Kazimierz \.Zorawski). It should also be noted that 1898  was a difficult year for Lie (who was in poor health). On May 22, 1898, he oficially resigned from his position in Leipzig in order to take up a special professorship in Kristiania (now Oslo). He was busy with writing up his research monographs and trying to have an input in naming his successor. In September 1898 Lie returned to Norway  and in February 1899 he died of pernicious anemia. (For information on the life and work of Sophus Lie and mathematics at Leipzig, cf. \cite{Fr}, \cite{Ko} and \cite{Stu}). Lucjan  B\"ottcher, Charles Bouton and Gerhard Kowalewski were his last doctoral students. \\

\section{B\"ottcher in Lw\'ow}

In 1901 B\"ottcher arrived in Lw\'ow and took a position of an \textit{asystent} in the Lw\'ow Polytechnic School .  His PhD diploma from Leipzig was nostrified (i.e., officially recognized) at the Lw\'ow  University. The same year (in a letter dated October 17) he applied for admission to habilitation at the Lw\'ow  University. Along with the diploma and curriculum vitae, he submitted offprints of two papers: ''Principles of iterational calculus, part three" and ''On properties of some functional determinants", as well as a plan  of lectures for 4 semesters. B\"ottcher's application was considered by a committee whose members were J\'ozef Puzyna\footnote{J\'ozef Puzyna (1856-1919) specialized in complex function theory and wrote the first textbook in Polish on the subject.}, Jan  Rajewski\footnote{Jan Rajewski (1857-1906) worked on differential equations and hypergeometric functions.}, Marian Smoluchowski\footnote{Marian Smoluchowski (1872-1917)  was an outstanding statistical physicist, working on Brownian motions  and diffusion.} and the dean Ludwik Finkel. The committee's decision, made on February 6, 1902, was not to admit B\"ottcher to habilitation. His scientific results were deemed correct but insufficient, although the committee also noted that the theory of iteration itself was not yet a well developed area of mathematics.  The recommendation was to wait until the official publication of  ''Principles of iterational calculus, part three" (or of potential works by B\"ottcher in areas of mathematics other than iteration and functional equations).\\
 
 In 1911 B\"ottcher obtained {\it veniam legendi} in mathematics in the c.k. Polytechnic School in Lw\'ow, where (since 1910) he was employed as an \textit{adiunkt}. In 1911 he  requested at the Faculty of Philosophy of the Lw\'ow University that his license to lecture at the c.k. Polytechnic School be also recognized at the university. His request was denied. \\
 
 Another time he applied for habilitation was in 1918. This time he had more publications to his name and he submitted six of them with the application. Here is the translation of (a fragment of) the committee's decision:\\
 
\textit{ ``The Candidate submitted along with the application the following works in mathematics:\\
  1) Major laws of convergence of iterations and their applications in analysis. Two papers in Russian, Kazan, 1903, 1905.\\
  2)  A note of solving the functional equation $\Psi f(z)-\Psi(z)=F(z)$, {\it Wiadomo\'sci Matematyczne}, vol. 13, Warsaw 1909\\
  3) Principles of iterational calculus, {\it Wektor}, 1912, Warsaw\\
  4) Nouvelle m\'ethode d'int\'egration d'un syst\`eme de $n$ \'equations fonctionelles lineair\'es du premier ordre de la forme $U_i(z)=\sum_{j=1}^{j=n}A_{i,j}(z)U_jF(z)$, Annales l'Ecole Normale Sup\'erieure, Paris, 1909\\
  5) A contribution to the calculus of iteration of a rational entire function, {\it Wiadomo\'sci Matematyczne}, vol.14, Warsaw 1912\\
  6) Iteration $f_x(z)$ of an algebraic function $f(z)$ metatranscendental in the index $x$, in Russian, Kazan 1912\\
``The paper no. 5 duplicates one written by the author in Polish and self-published already in 1905. (...) 
After formal deduction of formulas for a solution to the system of equations under investigation the author proceeds to give in \S  3 `A functional-theoretical discussion of the fundamental law', concluding boundedness of a certain set from finiteness of the numbers in it. Such reasoning is obviously erroneous, and therefore one cannot consider it to be proven that the series given by the author are-- under his conditions-- convergent.\\
``The method used by the Candidate in his works cannot be considered scientific. The author works with undefined, or ill-defined, notions (e.g., the notion of an iteration with an arbitrary exponent), and the majority of the results he achieves are transformations of one problem into another, no less difficult. In the proofs there are moreover illegitimate conclusions, or even fundamental mistakes. The author's popular, instructional works, e.g. `Principles of iterational calculus' ({\rm Wektor} 1912, no. 9, pp.  501-513, Warsaw), are written in an unclear manner. (...)\\ 
`` Despite great verve and determination, Dr. B\"ottcher's works do not yield any positive scientific results. There are many formal manipulations and computations in them; essential difficulties are usually dismissed with a few words without deeper treatment. The content and character diverges significantly from modern research.\\
``One should also add:\\
1. The shortcoming, or rather lack of rigor of the definition of iteration with an arbitrary exponent introduced by the candidate met with justified and clearly written criticism by Dr. Stanis\l aw Ruziewicz in {\rm Wektor}, Warsaw 1912, no. 5 [On a problem concerning commuting functions].\\
2. Dr. B\"ottcher applies for a second time for veniam legendi in mathematics. The first time the candidate was advised to withdraw his application because of the faults that the Committee at that time found with the candidate's works. These  faults and inadequacies were of the same nature which characterizes the candidate's work also today.\\
``The Committee's decision  passed unanimously on June 21, 1918: Not to admit Dr. B\"ottcher to further stages of habilitation. (signature illegible).}\\

B\"ottcher made his last attempt to obtain habilitation at the University  on May 1, 1919,  also unsuccessfully. He remained in the position of an \textit{adiunkt} at the Lw\'ow Polytechnic School until his retirement in 1935.\\


\section{The importance of B\"ottcher's work}

There are 19 known mathematical research publications by Lucjan Emil B\"ottcher. The following are the most important for the development of holomorphic dynamics:\\

\begin{enumerate}

\item Beitr\"age zu der Theorie der Iterationsrechnung, published by
Oswald Schmidt, Leipzig, pp.78, 1898 (doctoral dissertation).

\item Zasady rachunku iteracyjnego (cz\c e\'s\'c pierwsza i cz\c e\'s\'c druga)
[Principles of iterational calculus (part one and two)], Prace
Matematyczno Fizyczne, vol. X (1899 1900), pp. 65-86, 86-101

\item Zasady rachunku iteracyjnego (cz\c e\'s\'c III) [Principles of iterational
calculus (part III)], Prace Matematyczno Fizyczne, v.
XII(1901), p. 95-111

\item Zasady rachunku iteracyjnego (cz\c e\'s\'c III, doko\'nczenie) [Principles
of iterational calculus (part III, completion)], Prace Matematyczno Fizyczne, v. XIII(1902), pp. 353-371

\item Glavn"yshiye zakony skhodimosti iteratsiy i ikh prilozheniya k" 
analizu [The principal laws of convergence of iterates and their
application to analysis], Bulletin de la Societe Physico-Mathematique
de Kasan, tome XIII(1, 1903), p.137, XIV(2, 1904), p.
155-200, XIV(3, 1904), p. 201-234.

\end{enumerate}

In his doctoral thesis B\"ottcher set out to develop a general theory of iteration of functions with an arbitrary (not necessarily integer) exponent, in the newly available framework of Lie groups. He managed to formulate some basic properties and outline the relation with functional and differential equations, but almost halfway through he switched to the study of iteration of rational maps of the Riemann sphere. Unlike many of his predecessors working on iteration, he was interested in global rather than local behavior of maps. Here are his main ideas and results, expressed in modern terminology:\\

--the study of individual orbits of (iterated) rational maps, of their convergence and the limits that occur;\\
--the study of ``regions of convergence" (later called Fatou components) and their boundaries (Julia sets); determining the boundaries using backward iteration;\\
--relations between the local behavior of iterates and the magnitude of the derivative at a fixed point of a map;\\
--an example of an everywhere chaotic map, i.e., a map without regions of convergence constructed by means of elliptic functions;\\
--some observations about preperiodic points, which B\"ottcher called ``\.Zorawski points".\\

These topics re-emerged after 1918 in the works of other mathematicians and have since served as foundations of holomorphic dynamics. Fatou and Julia (independently) took advantage of the theory of normal families, formulated by Paul Montel, to   put forward and further study the division of the sphere (or the complex plane) into subsets in which iterates of a map display different behavior.  The limit maps of convergent subsequences of iterates and the behavior of their derivatives were studied in detail. An example of an everywhere chaotic map derived from an elliptic function was constructed  in 1918 by Samuel Latt\'es (hence the family of  maps  coming from such constructions is now known as Latt\'es examples).
The work of Kazimierz \.Zorawski,  preceding B\"ottcher's, has been  nearly forgotten, but  preperiodic points  nowadays enjoy renewed interest, due to their importance in the study of the parameter spaces of families of rational maps, including the famous Mandelbrot set (in relation to so-called Misiurewicz points, introduced by and named after Micha\l\ Misiurewicz, a Polish mathematician  active in the USA), as well as in arithmetic dynamics.\\

As for the relation between Lie's theory and iteration, the problem of defining iterations with an arbitrary exponent by embedding iterations of a  function into a one-parameter continuous group of transformation cannot be solved in such generality as B\"ottcher hoped for. For rational maps, I. N. Baker showed in 1960s that it is impossible to carry out such an embedding in the whole complex plane. Partially defined embeddings can be obtained in some cases.  On the other hand,  Julia also was interested in continuous groups of transformations, namely Kleinian groups. Many deep parallels between Kleinian groups and iterations of rational maps were observed and systematically explored in 1980s by Dennis Sullivan. \\

B\"ottcher himself called the combined papers (2), (3) and (4) a translation of his thesis, but they contain more 
results on (what is now known as) holomorphic dynamics than (1) and their organization is quite different. The new material, appearing mostly in (2), is the following:\\

--a different example of an everywhere chaotic map and a sketch of proof of its chaotic behavior;\\
--determining boundaries of the ``regions of convergence" for monomials and Chebyshev polynomials; study of simple dynamical properties of the ``boundary curves", e.g,, density of periodic points;\\
--brief attention given to irrationally neutral periodic points (without any conclusions);\\
--pointing out the role of critical points in the dynamics of rational maps; formulation of an exact upper bound for the number of (periodic) ``regions of convergence" in terms of the number of critical points of the map;\\ 
--the first formulation of B\"ottcher's theorem (about the local behavior of a map near a superattracting fixed point; cf. the introduction) and a sketch of its proof.\\

Again, the properties of ``boundary curves" (now known as Julia sets) were later studied in detail by Fatou and Julia (as well as by Salvatore Pincherle). Fatou also looked at irrationally neutral periodic points (pointing out that non-constant maps appear as limits of the sequences of iterates in a neighborhood of such a point) and examined the role of critical points in holomorphic dynamics. He  postulated the same upper bound on the number of ``regions of convergence" (more precisely, on the number of non-repelling periodic orbits) by the number of critical points as B\"ottcher did, 
 but he was able  to prove only a weaker one. The exact bound was finally proved by Mitsuhiro Shishikura in 1982; in 1999 Adam Epstein gave a different proof of the Fatou-Shishikura inequality.\\

The paper (5) is the most  cited publication by B\"ottcher. He again formulated his theorem in it and sketched its proof. Joseph Fels Ritt first cited this paper in  \cite{Ri}, where he wrote up a complete proof of B\"ottcher's theorem (unaware of its  earlier formulation in (2)), and other scholars have followed suit ever since.  The paper emphasizes relations between the theory of iteration and functional equations, and there are few new results in holomorphic dynamics in it, except a detailed analysis of the behavior of the map $z \mapsto z^2$  and the construction of an attracting basin by backward iteration. It is in this paper that the term ``chaotic" is introduced to describe the behavior of a map without regions of convergence. 
It should be noted that  ``Bulletin de la Societe Physico-Mathematique
de Kasan" was considered to be a prestigious journal, since Kazan was the home of  the International Lobachevsky Foundation awarding the Lobachevsky prize in geometry.\footnote{Sophus Lie was the first recipient of this prize in 1897 for his work on groups of transformations. In the years 1951-2000 the Lobachevsky prize continued to be awarded, first by the Soviet Academy of Sciences and since 1992 by the Russian Academy of Sciences. Kazan State University awards the Lobachevsky medal.} It  circulated widely, mainly through exchange for publications of other academic centers and learned societies; the volume X from year 1902 lists 123 institutions (51 in the Russian empire and 72 worldwide) participating in such exchange.\\

The papers (2), (3), (4) and (5), besides B\"ottcher's original contributions, contain detailed bibliography of related studies and an exhaustive discussion of  other matematicians' results, so B\"ottcher also comes out as very well versed in the literature of the subject. On the other hand, his papers contain many unjustified 
conclusions and some false statements, especially concerning iterates with an arbitrary exponent.  These shortcomings were noticed by the referees of his doctoral thesis in Leipzig and by the members of the habilitation committees in Lw\'ow.  But some of  inacurracies in B\"ottcher's paper are  typical  for an initial phase in the development of a new discipline, when the relevant  notions are just being formed.  E.g., he refers to the boundaries of regions of convergence as to ``boundary curves", while nowadays it is known that they can be totally disconnected--``dust-like" (such a situation occurs e.g. for a map $z \mapsto z^2+c$ with $c$ lying outside the Mandelbrot set).  However, in B\"ottcher's time set theory and topology were not advanced disciplines, and the notion of a curve was not clearly understood (Peano's example of a space-filling curve from 1890 was considered counterintuitive by some mathematicians). \\

B\"ottcher  revived his interest in iteration of rational maps around 1903 and subsequently published more  papers on this topic:\\

--Iteracye funkcyi liniowej [Iterations of a linear function], Wiadomo\'sci Matematyczne,
vol. VIII(1904), s. 291-307;\\
--Iteracye funkcyi liniowej (ci\c ag dalszy i doko\'nczenie) [Iterations of a linear function
(continuation and completion)], Wiadomo\'sci Matematyczne, vol. IX (1905), p. 77-86;\\
--Przyczynek do rachunku iteracyj funkcyi algebraicznej wymiernej ca\l kowitej [A contribution
to the calculus of iterations of an algebraic rational entire function], Wiadomo\'sci
Matematyczne XVI( 1912), s. 201-206;\\
--Iteracye funkcyi kwadratowej [Iterations of a quadratic function], Wiadomo\'sci Matematyczne,
 XVIII (1914), s. 83-132.\\

 Even though not always rigorous, B\"ottcher's mathematical output encompasses many ideas, examples and partial results that were later rediscovered independently by other mathematicians, giving rise to holomorphic dynamics as a new area of mathematics. One should hope that mathematicians and historians of science will recognize B\"ottcher's pioneering role in the formation of this discipline and will agree with the words of Alessandro Rosa:  ``Thus, at present, we have this `four of a kind' for global holomorphic dynamics: B\"ottcher, Fatou, Julia and Pincherle." (\cite{Ro}).\\

\textbf{Acknowledgments:} I thank Stanis\l aw Domoradzki for encouraging my involvement in the study of history of Polish mathematics and for constant sharing of ideas and materials related to Lucjan Emil B\"ottcher.  I also thank Alessandro Rosa for his remarks (in  lively email exchanges) on the contents of B\"ottcher's work and its significance.  Finally, I thank Terry Czubko for her help with reading handwritten German documents.

\textbf{Dedication:} I dedicate this article to  my grade school mathematics teachers: \\ Maria Burek, Janina \'Sl\'osarczyk, Maria Kubin and Ewa Dutkiewicz.

\end{document}